\documentclass[11pt]{article}
\usepackage{amssymb}

\newcommand{\Z}{{\mathbb Z}}
\newcommand{\C}{{\mathbb C}}

\newcommand{\F}{{\mathbb F}}
\newcommand{\eps}{\varepsilon}
\newcommand{\EE}{{\mathbb E}}

\newcommand{\SIGMA}{\raisebox{-0.4ex}{\mbox{\Large $\Sigma$}}}

\newtheorem{theorem}{Theorem}
\newtheorem{lemma}{Lemma}

\title{On sumsets and spectral gaps}

\author{Ernie Croot\thanks{Supported in part by an NSF grant.} \\ Georgia Tech \\
School of Mathematics \\ 103 Skiles \\ Atlanta, GA 30332 \\ \\
Tomasz Schoen\thanks{Research partially supported by MNSW grant 2
P03A 029 30}
\\ Department of Discrete Mathematics \\ Adam Michiewicz University \\
ul. Umultowska 87, 61-614 Pozna\'n, Poland}

\begin{document}

\maketitle

\begin{abstract} Suppose that $S \subseteq \F_p$, where $p$ is a prime number.
Let $\lambda_1,...,\lambda_p$ be the Fourier coefficients of $S$ arranged as follows
$$
|\hat S(0)|\ =\ |\lambda_1|\ \geq\ |\lambda_2|\ \geq\ \cdots\ \geq\ |\lambda_p|.
$$
Then, as is well known, the smaller $|\lambda_2|$ is, relative to $|\lambda_1|$,
the larger the sumset $S+S$ must be; and, one can work out as a function of
$\eps$ and the density $\theta = |S|/p$,
an upper bound for the ratio $|\lambda_2|/|\lambda_1|$ needed in order
to guarantee that
$S+S$ covers at least $(1-\eps)p$ residue classes modulo $p$.   Put another way,
if $S$ has a large spectral gap, then most elements of $\F_p$
have the same number of representations as a sum of two elements of $S$,
thereby making $S+S$ large.

What we show in this paper is an extension of this fact, which holds for
spectral gaps between other consecutive Fourier coefficients
$\lambda_k,\lambda_{k+1}$, so long as $k$ is not too large; in particular, our
theorem will work so long  as
$$
1\ \leq\ k\ <\ {\log p \over \log 4}
$$
Furthermore, we develop results for repeated sums $S+S+\cdots +S$.

It is worth noting that this phenomena does not hold in the larger
finite field setting $\F_{p^n}$ for fixed $p$, and where we let $n
\to \infty$, because, for example, the indicator function for a
large subspace of $\F_{p^n}$ can have a large spectral gap, and yet
the sumset of that subspace with itself equals the subspace (which
therefore means it cannot cover density $1-\eps$ fraction of
$\F_{p^n}$). The property of $\F_p$ that we exploit, which does not
hold for $\F_{p^n}$ (at least not in the way that we would like --
Browkin, Divis and Schinzel \cite{browkin} have analyzed the problem
for more general settings than just $\F_p$), is something we call a
``unique differences'' property, first identified by W. Feit, with
first proofs and basic results found by Straus \cite{straus}.
\end{abstract}

\section{Introduction}

Supose that
$$
f\ :\ \F_p\ \to\ [0,1],
$$
and let
$$
\theta\ :=\ \EE(f)\ :=\ p^{-1} \SIGMA_n f(n).
$$
For an $a \in \F_p$, define the usual Fourier transform
$$
\hat f(a)\ :=\ \SIGMA_n f(n) e^{2\pi i an/p}.
$$
We order the elements of $\F_p$ as
$$
a_1,...,a_p,
$$
so that
\begin{equation} \label{ai_def}
|\hat f(a_1)|\ \geq\ |\hat f(a_2)|\ \geq\ \cdots\ \geq\ |\hat f(a_p)|;
\end{equation}
(there may be multiple choices for $a_1,...,a_p$ -- {\it any} ordering will do)
and, for convenience, we set
$$
\lambda_i\ =\ f(a_i),\ i=1,...,p.
$$
Note, then, that
$$
|\lambda_1|\ =\ |\hat f(0)|.
$$
\bigskip

In this paper we prove the following basic theorem.

\begin{theorem} \label{main_theorem}  Suppose that
$f : \F_p \to [0,1]$, $f$ not identically $0$, has the property that for some
$$
1 \leq k < {\log p \over \log 4}
$$
we have that
$$
|\lambda_{k+1}|\ \leq\ \gamma |\lambda_k|.
$$
Then,
$$
|\{ n \in \F_p\ :\ (f*f)(n)\ >\ 0\}|\ \geq\ p (1 -  2\theta p^2 \gamma^2 |\lambda_k|^{-2}).
$$
\end{theorem}

\noindent {\bf Remark.}  It is easy to construct functions $f$ which have a large
spectral gap as in the hypotheses.  For example, take $f$ to be the function
whose Fourier transform satisfies $\hat f(0) = p/2$ and $\hat f(1) = \hat f(-1) = p/4$,
then $\hat f(a) = 0$ for $a \neq 0, \pm 1$.  Clearly we have $f : \F_p \to [0,1]$,
and of course $f$ has a large spectral gap between $\lambda_3$ and $\lambda_4$
($|\lambda_3| = p/4$, while $\lambda_4 = 0$).
\bigskip

By considering repeated sums, one can prove similar sorts of results, but
which hold for a much wider range of $k$.  Furthermore, one can derive conditions
guaranteeing that $(f*f*\cdots *f)(n) > 0$ for all $n \in \F_p$, not just $1-\eps$
proportion of $\F_p$; and, these conditions are much simpler and cleaner than
those of Theorem \ref{main_theorem} above.   This new theorem is given as follows:

\begin{theorem} \label{main_theorem2}  Suppose that $f : \F_p \to [0,1]$,
$f$ not identically $0$, has the property that for some
$$
1\ \leq\ k\ <\ (\log p)^{t-1} (5t \log\log p)^{-2t+2},
$$
we have that
$$
|\lambda_{k+1}|\ \leq\ \gamma |\lambda_k|,\ {\rm where\ }
\gamma\ <\ t^{-1} \theta^{-t+2} (|\lambda_k|/p)^{t-1}.
$$
Then, for $t \geq 3$, the $t$-fold convolution $f*f*\cdots *f$ is positive on all of $\F_p$.
\end{theorem}

\noindent {\bf Remark.}  It is possible to prove even stronger results for when $k$
is much smaller than $t$ (say less than the square-root of $t$),
though the result is a little more technical to state.
\bigskip

We conjecture that it is possible to prove a lot more:
\bigskip

\noindent {\bf Conjecture.}  It is possible to develop bounds of the same general
quality as to those in Theorem \ref{main_theorem} for the number of
$n$ with $(f*f)(n) > 0$, given that $f$ has a large spectral gap between the
$k$th and $(k+1)$st largest Fourier coefficients of $f$, for any
$k < p^{1/2}$, say.   This would obviously require a different sort of proof than
appears in the present paper, as a key lemma we use
(Lemma \ref{unique_differences_lemma}) is close to best-possible.
Furthermore, it should be possible to prove a version of
Theorem \ref{main_theorem2} under the assumption of such a spectral gap.

\section{Some lemmas}

\begin{lemma}[Dirichlet's Box Principle] \label{dirichlet}  Suppose that
$$
r_1,...,r_t\ \in\ \F_p.
$$
Then, there exists non-zero $m \in \F_p$ such that
$$
{\rm For\ } i =1,...,t,\ \left | \left | {m r_i \over p} \right | \right |\ \leq\ p^{-1/t},
$$
where here $||x||$ denotes the distance from $x$ to the nearest integer.
\end{lemma}

The proof of this lemma is standard, so we omit it.  The following
lemma is also standard, and was first discovered by Straus
\cite{straus} (and re-discovered by the first author) though we will
bother to give the proof.  It is worth remarking that Browkin, Divis
and Schinzel \cite{browkin} have worked out a more general version
of this lemma that holds in artibrary groups; and,  Lev \cite{lev}
has extended and applied these results to address some problems on
discrepancy.

\begin{lemma}[Unique Differences Lemma]  \label{unique_differences_lemma}
Suppose that
$$
B\ :=\ \{b_1,...,b_t\} \subseteq\ \F_p.
$$
Then, if
$$
t\ <\ (\log p)/\log 4,
$$
there will exist $d \in \F_p$ having a unique representation as a difference of
two elements of $B$.
\end{lemma}

\noindent {\bf Proof of the lemma.}  First, from the Dirichlet Box Principle above,
we deduce that there exists a non-zero dilation constant $m \in \F_p$ such that
if we let
$$
c_i\ \equiv\ m b_i \pmod{p},\ |c_i|\ <\ p/2,
$$
then, in fact,
$$
|c_i|\ \leq\ p^{1-1/t}.
$$
So long as
$$
p^{1-1/t}\ <\ p/4\ \ \iff\ \ p\ >\ 4^t,
$$
we will have that all these $c_i$ lie in $(-p/4, p/4)$.  Then, if we let
$$
c_x\ :=\ \min_i c_i,\ {\rm and\ } c_y\ :=\ \max_i c_i,
$$
we claim that $d \in B-B$ given by
$$
d\ =\ c_y - c_x
$$
has a unique representation as a difference of elements of $B$, and therefore
$c_y - c_x$ is that unique representation.   The reason that this is the case is that
since $c_i \in (-p/4,p/4)$ we have that all the differences
$$
c_i - c_j\ \in\ (-p/2, p/2);
$$
and so, two of these differences are equal if and only if they are equal modulo $p$;
and, it is clear that, over the integers, $d = c_y - c_x$ has a unique representation,
implying that it has a unique representation modulo $p$.
\hfill $\blacksquare$
\bigskip

We will actually need a generalization of this lemma, which is a refinement of
one appearing in \cite{luczak}, and is given as follows.

\begin{lemma}  \label{small_rep} Suppose that
$$
B_1,B_2\ \subseteq\ \F_p,
$$
where
$$
10\ \leq\ |B_1|\ \leq\ p/2,\ {\rm and\ } 3 |B_2| \log |B_1| > \log p.
$$
Then, there exists $d \in B_1-B_2$ having at most
$$
20|B_2| (\log |B_1|)^2/\log p
$$
representations as
$$
d\ =\ b - b',\ b \in B_1, b' \in B_2.
$$

Furthermore, if
$$
1\ \leq\ |B_1|\ \leq\ p/2,\ {\rm and\ } 3 |B_2| \log |B_1|\ <\ \log p,
$$
then there exists $d \in B_1 - B_2$ having a unique representation as
$d = b_1 - b_2$, $b_1 \in B_1, b_2 \in B_2$.
\end{lemma}

\noindent {\bf Proof of the lemma.}  Let $B'$ be a random subset of $B_2$, where
each element $b \in B_2$ lies in $B'$ with probability
$$
(\log p)/3 |B_2| \log |B_1|.
$$
Note that this is where our lower bound $3 |B_2| \log |B_1| > \log p$ comes in,
as we need this to be less than $1$.

So long as the $B'$ we choose satisfies
\begin{equation} \label{b'_bound}
|B'|\ <\ (\log p)/2 \log |B_1|,
\end{equation}
which it will with probability at least $1/2$,
we claim that there will always exist an element
$d \in B-B'$ having a unique representation as a difference $b_1-b_2'$,
$b_1 \in B, b_2' \in B'$:  First, note that it suffices to prove this for the
set $C_1-C'$, where
$$
C_1 = m\cdot B_1,\ C_2 = m\cdot B_2,\ {\rm and\ } C' = m\cdot B',
$$
where $m$ is a dilation constant
chosen according to Dirichlet's Box Lemma so that every element $x \in C'$
(when considerecd as a subset of $(-p/2, p/2]$) satisfies
$$
|x|\ \leq\ p^{1- 1/|B'|}\ <\ p/3|B_1|.
$$
Now, there must exist an integer interval
$$
I\ :=\ (u,\ v)\ \cap\ \Z,\ u,v \in C_1,
$$
(which we consider as an interval modulo $p$) such that
$$
|I|\ \geq\ p/|C_1| - 1\ =\ p/|B_1| - 1,
$$
and such that no element of $C_1$ is congruent modulo $p$ to an element of $I$.
Clearly, then, one of the following two elements
$$
v - \max_{c' \in C'} c',\ {\rm or\ } u - \min_{c' \in C'} c'
$$
(here, this $c'$ is thought of an an element of $(-p/2, p/2]$)
has a unique representation as a difference.  The reason we need this either-or
is that all the elements of $C'$ could be negative.

Now we define the functions
\begin{eqnarray*}
\nu(x)\ &:=&\ |\{ (c_1, c_2) \in C_1 \times C_2\ :\ c_1 - c_2\ =\ x\}|;\ {\rm and,} \\
\nu'(x)\ &:=&\ |\{ (c_1,c_2') \in C_1 \times C'\ :\ c_1 - c_2'\ =\ x\}|.
\end{eqnarray*}
We claim that with probability at least $1/2$ we will have that
\begin{eqnarray} \label{implication}
{\rm for\ every\ }x\in \F_p,\ \nu(x) > 20 |B_2| (\log |B_1|)^2/\log p\
\Longrightarrow\ \nu'(x) \geq 2. \nonumber \\
\end{eqnarray}
To see this, fix $x \in C_1-C_2$.  Then, $\nu'(x)$ is the following
sum of independent Bernoulli random variables:
$$
\nu'(x)\ =\ \sum_{j=1}^{\nu(x)} X_j,\ {\rm where\ } {\rm Prob}(X_j = 1)\ =\
(\log p)/3 |B_2| \log |B_1|.
$$
The variance of $\nu'(x)$ is
$$
\sigma^2\ =\ \nu(x) {\rm Var}(X_1)\ \leq\ \nu(x) \EE(X_1).
$$
We now will need the following well-known theorem of Chernoff:

\begin{theorem} [Chernoff's inequality]  Suppose that
$Z_1,...,Z_n$ are independent
random variables such that $\EE(Z_i) = 0$ and $|Z_i| \leq 1$ for all $i$.  Let
$Z := \SIGMA_i Z_i$, and let $\sigma^2$ be the variance of $Z$.  Then,
$$
{\rm Prob}(|Z| \geq \delta \sigma)\ \leq\ 2 e^{-\delta^2/4},\ {\rm for\ any\ }
0 \leq \delta \leq 2\sigma.
$$
\end{theorem}

We apply this theorem using $Z_i = X_i - \EE(X_i)$ and
$$
\delta\sigma = \nu(x) \EE(X_1) - 1.
$$
and then quickly deduce that if $\nu(x) > 20 |B_2| (\log |B_1|)^2/\log p$, then
$$
{\rm Prob}(\nu'(x) \leq 1)\ \leq\ 2 \exp \left ( -\delta^2/4\right )\ <\ 1/2|B_1|,
$$
for $p$ sufficiently large.  Clearly, then, with probability at least $1/2$
we will have that (\ref{implication}) holds for all $x$, as claimed.  But we also
had that (\ref{b'_bound}) holds with probability at least $1/2$; so,
there is an instantiation of the set $B'$ such that {\it both} (\ref{implication})
and (\ref{b'_bound}) hold.  Since we proved that such
$B'$ has the property that there is an
element of $x \in B_1-B'$ having $\nu'(x) = 1$, it follows from (\ref{implication}) that
$\nu(x) \leq 20|B_2| (\log |B_1|)^2/\log p$, which proves the first part of our
lemma.

Now we prove the second part of the lemma:  First, the lemma is
obviously true in the case $|B_1| = 1$, so we assume that $|B_1|
\geq 2$. Since we are also assuming that $|B_2| < \log p/3 \log
|B_1|$, we have by the Dirichlet Box Principle there exists $m$ such
that for every $x \in C_2 := m\cdot B_2$ we have $|x| \leq
p/|B_1|^3$; furthermore, by the pigeonhole principle there exists an
integer interval $I := (u,v) \cap \Z$ with $u,v \in C_1 := m\cdot
B_1$, with $|I| \geq p/|B_1| - 1$, which contains no elements of
$B_1$. So, either
$$
v - \max_{x \in C_2} x\ {\rm or\ } u - \min_{x \in C_2} x
$$
has a unique representation as a difference $c_1-c_2$, $c_1 \in C_1$, $c_2 \in C_2$.
The same holds for $B_1 - B_2$, and so our lemma is proved.

\hfill $\blacksquare$

\section{Proof of Theorem \ref{main_theorem}}

We apply this last lemma with
$$
B\ =\ A\ =\ \{a_1,...,a_k\},\ {\rm so\ } t =k.
$$
Then, let $d$ be as in the lemma, and let
$$
a_x,\ a_y\ \in\ A
$$
satisfy
$$
a_y - a_x\ =\ d.
$$

We define
$$
g(n)\ :=\ e^{2\pi i dn/p} f(n),
$$
and note that
$$
(f*f)(n)\ \geq\ |(g*f)(n)|
$$
So, our theorem is proved if we can show that $(g*f)(n)$ is often non-zero.
Proceeding in this vein, let us compute the Fourier transform of
$g*f$:  First, we have that
$$
\hat g(a)\ =\ \SIGMA_n g(n) e^{2\pi i an/p}\ =\ \SIGMA_n f(n) e^{2\pi i n(a+d)/p}
\ =\ \hat f(a+d).
$$
So, by Fourier inversion,
\begin{equation} \label{fg}
(f*g)(n)\ =\  p^{-1} e^{-2\pi i a_x/p} \hat f(a_x) \hat f(a_y)\ +\  E(n),
\end{equation}
where $E(n)$ is the ``error'' given by
$$
E(n)\ =\ p^{-1} \SIGMA_{i \neq x} e^{-2\pi i a_i n/p} \hat f(a_i) \hat f(a_i + d).
$$
Note that for every value of $i \neq x$ we have that
\begin{eqnarray} \label{eitheror}
&& {\rm either\ } a\ \ {\rm or\ \ } a + d\ \ {\rm lies\ \ in\ \ }
\{a_{k+1},...,a_p\} \nonumber \\
&&\hskip0.5in \Longrightarrow\ \ |\hat f(a)\hat f(a + d)|\ \leq\
\gamma |\lambda_k| \max \{ |\hat f(a)|,\ |\hat f(a+d)|\}. \nonumber \\
\end{eqnarray}

To finish our proof we must show that ``most of the time'' $|E(n)|$ is smaller than
the ``main term'' of (\ref{fg}); that is,
$$
|E(n)|\ <\ p^{-1} |\hat f(a_x) \hat f(a_y)|.
$$
Note that this holds whenever
\begin{equation} \label{fg_target}
|E(n)|\ <\ p^{-1} |\lambda_k|^2.
\end{equation}

We have by Parseval and (\ref{eitheror}) that
\begin{eqnarray*}
\SIGMA_n |E(n)|^2\ &=&\ p^{-1} \SIGMA_{i \neq x} |\hat f(a_i)|^2 |\hat f(a_i + d)|^2 \\
&\leq&\ 2p^{-1} \gamma^2 |\lambda_k|^2 \SIGMA_a |\hat f(a_i)|^2 \\
&\leq&\ 2\gamma^2 |\lambda_k|^2 \hat f(0).
\end{eqnarray*}
So, the number of $n$ for which (\ref{fg_target}) holds is at least
$$
p(1\ -\ 2\gamma^2 |\lambda_k|^{-2} \hat f(0) p)\ =\
p (1 -  2p^2 \theta \gamma^2 |\lambda_k|^{-2}),
$$
as claimed.

\section{Proof of Theorem \ref{main_theorem2}}

Let
$$
B_1\ :=\ B_2\ :=\ A\ =\ \{a_1,...,a_k\}.
$$

Suppose initially that $3|A| \log |A| > \log p$, so that the hypotheses of the
first part of Lemma \ref{small_rep} hold.  We have then that there
exits $d_1\in B_1-B_2 = A-A$ with at most
$20 |A| (\log |A|)^2/\log p$ representations as
$d_1 = a-b$, $a,b \in A$.   Let  now $A_1$ denote the set of all the
elements $b$ that occur.  Clearly,
$$
|A_1|\ \leq\ 20 |A| (\log |A|)^2/\log p.
$$

Keeping $B_1 = A$, we reassign $B_2 = A_1$.   So long as
$3 |A_1|\log |A| > \log p$ we may apply the first part of
Lemma \ref{small_rep}, and when we
do we deduce that there exists $d_2 \in A-A_1$ having at most
$20 |A_1| (\log |A|)^2/\log p$ representations as $d_2 = a-b$, $a\in A$, $b \in A_1$.
Let now $A_2$ denote the set of all elements $b$ that occur.  Clearly
$$
|A_2|\ \leq\ 20 |A_1| (\log |A|)^2/\log p.
$$

We repeat this process, reassigning $B_2 = A_2$, then $B_2 = A_3$, and so on,
all the while producing these sets $A_1,A_2,...$ and differences
$d_1,d_2,...$, until we reach a set $A_m$ satisfying
$$
3 |A_m| \log |A|\ <\ \log p.
$$
We may, in fact, reach this set $A_m$ with $m=1$ if $3|A| \log |A| < \log p$.

It is clear that since at each step we have
$$
|A_i|\ \leq\ 20 |A_{i-1}| (\log |A|)^2/\log p,
$$
and since we have assumed that
$$
|A|\ <\ (\log p)^{t-1} (5t \log\log p)^{-2t+2},
$$
we will reach such a set with $m$ of size at most
$$
m\ \leq\ t-1.
$$

This set $A_m$ will have the property, by the second part of Lemma \ref{small_rep},
that there exists $d_m \in A-A_m$ having a unique representation as $d_m = a-b$,
$a \in A$, $b \in A_m$.

Now, we claim that there exists unique $b \in \F_p$ such that
$$
b,\ b + d_1,\ b + d_2,\ ...,\ b+d_m\ \in\ A.
$$
To see this, first let $b \in A$.  Since $b + d_1 \in A$ we must
have that $b \in A_1$, by definition of $A_1$.  Then, since $b + d_2 \in A$, it follows
that $b \in A_2$.  And, repeating this process, we eventually conclude that
$b \in A_m$.

So, since $b \in A_m$, and $b + d_m \in A$, we have $d_m = a-b$, $a \in A$, $b \in A_m$.
But this $d_m$ was chosen by the second part of Lemma \ref{small_rep} so that
it has a unique representation of this form.  It follows that $b \in A$ is unique, as
claimed.
\bigskip

From our funciton $f : \F_p \to [0,1]$, we define the functions $g_1,g_2,...,g_m : \F_p \to \C$
via
$$
f_i(n)\ :=\ e^{2\pi i d_i n/p} f(n).
$$
It is obvious that
$$
{\rm support}(f*f*\cdots f *g_1*g_2*\cdots *g_m)\ \subseteq\ {\rm support}(f*f*\cdots * f),
$$
where there are $t$ convolutions on the left, and $t$ on the right; so, $f$ appears
$t-m$ times on the left.

We also have that
$$
\hat g_i(a)\ =\ \hat f(a+d_i),
$$
and therefore
$$
(\widehat {f*f*\cdots * f*g_1*\cdots*g_m})(a)\ =\
\hat f(a)^{t-m} \hat f(a+d_1) \hat f(a+d_2) \cdots \hat f(a+d_m).
$$
Since there exists unique $a$, call it $x$, such that all these $a+d_i$ belong to $A$,
we deduce via Fourier inversion that for any $n \in \F_p$,
$$
(f*f*\cdots *g_1*\cdots * g_m)(n)\ =\ p^{-1} e^{-2\pi i nx/p} \hat f(x)^{t-m}
\hat f(x+d_1) \cdots \hat f(x+d_m)\ +\ E(n),
$$
where the ``error'' $E(n)$ satisfies, by the usual $L^2-L^\infty$ bound,
$$
|E(n)|\ \leq\ t |\lambda_{k+1}| \theta^{t-3} p^{t-4} \SIGMA_a |\hat f(a)|^2\ \leq\
t \gamma (\theta p)^{t-2} |\lambda_k|.
$$
So, whenever this is smaller than that main term, we have that the convolution is
non-zero, and therefore so is $(f*f*\cdots *f)(n)$.  This occurs if
$$
t \gamma (\theta p)^{t-2} |\lambda_k|\ \leq\ p^{-1} |\lambda_k|^t,
$$
which holds whenever $t \geq 2$ and
$$
\gamma\ <\ t^{-1} \theta^{-t+2} (|\lambda_k|/p)^{t-1}.
$$


\begin{thebibliography}{999}

\bibitem{browkin} J. Browkin, B. Divis, and A. Schinzel,
{\it Addition of sequences in general fields}, Monatsh. Math.
{\bf 82} (1976), 261-268.

\bibitem{lev} V. Lev, {\it Simulatenous approximations and covering by
arithmetic progressions in $\F_p$}, Jour. Comb. Theory Ser. A {\bf 92} (2000),
103-118.

\bibitem{luczak} T. {\L uczak} and T. Schoen, {\it On a problem of Konyagin}, preprint.

\bibitem{straus} E. G. Straus, {\it Differences of residues ${\rm mod} p$},
J. Number Theory {\bf 8} (1976), 40-42.


\end{thebibliography}
\end{document}